 \documentclass[12pt,twoside]{amsart}
 \usepackage{amssymb}
 \usepackage{amscd}

 \title[Equivariant completions]
 {Equivariant completions of toric contraction morphisms}  
 \author{Osamu Fujino} 
 \subjclass[2000]{Primary 14M25; Secondary 14E30.}
 \date{2005/3/26}
 \keywords{Toric varieties, Mori theory, Minimal Model Program, 
 Equivariant completion}
\address{Graduate School of Mathematics\\ 
Nagoya University, Chikusa-ku Nagoya 464-8602 Japan}
 \email{fujino@math.nagoya-u.ac.jp}
\newcommand{\Pic}[0]{{\operatorname{Pic}}}
\newcommand{\Proj}[0]{{\operatorname{Proj}}}
\newcommand{\red}[0]{{\operatorname{red}}}
\newtheorem{thm}{Theorem}[section]
\newtheorem{lem}[thm]{Lemma}
\newtheorem{cor}[thm]{Corollary}
\newtheorem{prop}[thm]{Proposition}
\newtheorem*{claim}{Claim}

\theoremstyle{definition}
\newtheorem{ex}[thm]{Example}
\newtheorem{defn}[thm]{Definition}
\newtheorem{que}[thm]{Question}
\newtheorem{rem}[thm]{Remark}
\newtheorem*{ack}{Acknowledgments}       
 
\newtheorem{say}[thm]{}
\newtheorem*{notation}{Notation}         
       
\begin{document}
\bibliographystyle{amsalpha+}

\begin{abstract}
We treat equivariant completions of 
toric contraction morphisms 
as an application of the toric Mori theory. 
For this purpose, 
we generalize the toric Mori theory for 
non-$\mathbb Q$-factorial toric varieties.  
So, our theory seems to be quite different from 
Reid's original combinatorial toric Mori theory. 
We also explain various examples of 
non-$\mathbb Q$-factorial 
contractions, which imply that the $\mathbb Q$-factoriality 
plays an important role in the Minimal Model Program. 
Thus, this paper completes the foundation of 
the toric Mori theory and shows us a new aspect of 
the Minimal Model Program. 
\end{abstract}

\maketitle

\section{Introduction}\label{intro}
In \cite{part1}, we gave a simple and 
non-combinatorial proof 
to the toric Mori theory. As mentioned in \cite{part1}, our 
method cannot recover combinatorial 
aspects of \cite{reid}. 
One of the 
main purposes of this paper is to understand the local 
behavior of the toric contraction morphisms, 
which was described in \cite[(2.5) Corollary]{reid} when 
the varieties are $\mathbb Q$-factorial and 
{\em{complete}}. 
It is obvious that the non-complete fans are much harder 
to treat than the complete ones. 
So, we avoid manipulating and subdividing 
non-complete fans. 
Our strategy is to 
compactify the toric contraction morphisms equivariantly and 
apply Reid's result. 

Let $f:X\longrightarrow Y$ be a projective toric morphism. 
We would like to compactify $f:X\longrightarrow Y$ equivariantly, 
that is, 
$$
\begin{matrix}
\overline f: & \overline X& \longrightarrow &\overline Y\\
                   & \cup & &\cup\\
f: & X& \longrightarrow &Y,\\
\end{matrix}
$$
where $\overline X$ (resp.~$\overline Y$) 
is an equivariant completion of $X$ (resp.~$Y$). 
More precisely, 
we would like to compactify $f$ equivariantly without 
losing the following properties: 
\begin{enumerate}
\item[(i)] projectivity of the morphism, 
\item[(ii)] $\mathbb Q$-factoriality of the source space, 
\item[(iii)] the relative Picard number is one, 
\end{enumerate}
and so on. Note that we do not assume 
$f$ to be birational. 
The main results are 
Theorems \ref{kore} and \ref{supple-ne}, where we compactify 
$f$ equivariantly by using the toric Mori theory. 
The statements are too long to mention here. 
These theorems guarantee that 
we can always compactify toric contraction morphisms 
equivariantly preserving nice properties. 
However, our proof does not show us how to 
compactify $f$ even if it is given 
concretely. 
As a corollary, we obtain a description of 
the toric contraction morphisms when 
the source spaces are $\mathbb Q$-factorial and 
the relative Picard numbers are one (Theorem \ref{desc}). 
As mentioned above, 
it seems to be difficult to 
obtain a local description of the toric contraction 
morphism without reducing it to the complete case. 
This is why we treat the equivariant completions of 
the toric contraction morphisms. 

To carry out our program, 
we generalize the toric Mori theory for 
{\em{non-complete}} and 
{\em{non-$\mathbb Q$-factorial}} varieties. 
It is also the main theme of this paper. 
We do not need the toric Mori theory for 
non-$\mathbb Q$-factorial varieties to construct 
an equivariant completion of $f:X\longrightarrow Y$ 
when $X$ is $\mathbb Q$-factorial. However, 
it is natural and interesting to consider 
the toric Mori theory for non-$\mathbb Q$-factorial 
varieties. 
Without $\mathbb Q$-factoriality, various new phenomena 
occur even in the three-dimensional Minimal Model Program 
(see Section \ref{sa} and \cite{fs-flop2}). We believe that this generalized 
version of the toric Mori theory is not reachable by Reid's 
combinatorial technique since non-complete and non-simplicial 
fans are very difficult to manipulate. 
So, it seems to be reasonable to regard our toric Mori theory 
to be different from Reid's combinatorial one. 
The coverage of our theory is much wider than Reid's. 

Note that 
the Minimal Model Program for non-$\mathbb Q$-factorial 
varieties may be useful in the study of higher dimensional 
log flips (see \cite[Section 4]{fu-pri}). 
This paper will open the door to the non-$\mathbb Q$-factorial 
world. 

This paper mainly treats the conceptual 
aspects of the toric Mori theory. 
The appendix \cite{fs-flop2} (see also \cite{fs-flop}), 
where we construct an 
example of global toric $3$-dimensional flops, will 
supplement this paper from the combinatorial viewpoint. 
Recently, Hiroshi Sato described the combinatorial 
aspects of the toric contraction morphisms from 
$\mathbb Q$-factorial toric varieties by 
using {\em{extremal primitive relations}} and 
Theorem \ref{kore}. For the details, see \cite{sato}. 
It will help us to understand Theorem \ref{desc} below. 

\vspace{3mm}

We summarize the contents of this paper:~In 
Section \ref{sec2}, 
we prove the existence of equivariant completions of toric contraction morphisms in various settings. 
For this purpose, 
we generalize the toric Mori 
theory for non-$\mathbb Q$-factorial toric varieties. 
Section \ref{sec-app} deals with applications of the 
equivariant completions obtained in Section \ref{sec2}. 
The final theorem in Section \ref{sec-app} 
is a slight generalization of 
the main theorem of \cite{fujino}. 
In Section \ref{sa}, we will treat various examples of 
non-$\mathbb Q$-factorial 
toric contraction morphisms. They imply that 
it is difficult to describe the local behavior of the (toric) 
contraction morphisms without the $\mathbb Q$-factoriality 
assumption. 
This section is independent of the other sections and 
seems to be valuable for those studying the Minimal 
Model Program.  

\begin{ack}
The author would like to thank Professors J\'anos Koll\'ar, 
Masanori Ishida for comments, and Florin Ambro for 
pointing out a mistake. 
He also likes to thank Dr.~Hiroshi Sato for constructing 
a beautiful example and pointing out some mistakes. 
Thanks are due to the Institute for Advanced Study 
for hospitality. 
He was partially supported by a grant from the 
National Science Foundation:~DMS-0111298. 
Finally, the author thanks the referee and Professor 
Kenji Matsuki, whose comments helped him 
to correct errors. 
\end{ack}

\begin{notation}
We often use the notation and the results in 
\cite{part1}. 
We will work over an algebraically closed field $k$ 
throughout this paper. 

(i) Let $v_i\in N\simeq \mathbb Z^n$ for 
$1\leq i\leq k$. Then the symbol $\langle v_1,v_2, \cdots, v_k\rangle$ 
denotes the cone $\mathbb R_{\geq 0}v_1+\mathbb R_{\geq 0}v_2+\cdots 
+\mathbb R_{\geq 0}v_k$ in $N_{\mathbb R}\simeq \mathbb R^n$, 
where $\mathbb R_{\geq 0}$ is the set of non-negative real numbers.

(ii) A {\em{toric morphism}} $f:X\longrightarrow Y$ 
means an equivariant morphism $f$ between toric varieties 
$X$ and $Y$.   
\end{notation}

\section{Equivariant completions of toric contraction 
morphisms}\label{sec2}

Let us start with the following preliminary proposition. 
Its proof is a warm-up of our toric Mori theory \cite{part1}. 

\begin{prop}\label{sai}
Let $f:X\longrightarrow Y$ be a projective 
toric morphism and $\overline Y$ an equivariant completion of $Y$. 
Then there 
exists an equivariant completion 
of $f:X\longrightarrow Y${\em{;}} 
$$
\begin{matrix}
\overline f: & \overline X& \longrightarrow &\overline Y\\
                   & \cup & &\cup\\
f: & X& \longrightarrow &Y,\\
\end{matrix}
$$
where 
\begin{enumerate}
\item[(i)] $\overline X$ is an equivariant completion of $X$, 
and 
\item[(ii)] $\overline f$ is a {\em{projective}} toric morphism. 
\end{enumerate}
Furthermore, 
\begin{enumerate}
\item[(1)] if $X$ is $\mathbb Q$-factorial $($see \cite
[Definition 2.3]{part1}$)$, 
then we can make $\overline X$ to be 
$\mathbb Q$-factorial, and 
\item[(2)] if $X$ has only $($$\mathbb Q$-factorial$)$ 
terminal $($resp.~canonical$)$ 
singularities $($see \cite[(1.11) Definition]{reid} or 
\cite[Definition 2.9]{part1}$)$, 
then we can make $\overline 
X$ to have only $($$\mathbb Q$-factorial$)$ 
terminal $($resp.~canonical$)$ singularities.  
\end{enumerate}
\end{prop}
\begin{proof}
By Sumihiro's equivariant embedding theorem, 
there exists an equivariant completion $X_1$ 
of $X$. 
Let $X_2$ be the graph of the rational map $f:X_1
\dashrightarrow 
\overline Y$. Then, we obtain 
$$
\begin{matrix}
f_2: & X_2& \longrightarrow &\overline Y\\
                   & \cup & &\cup\\
f: & X& \longrightarrow &Y. \\
\end{matrix}
$$
Let $D$ be an $f$-ample Cartier divisor on $X$ and $D_2$ 
the closure of $D$ on $X_2$. 
By Corollary 5.8 
in \cite{part1}, 
$\bigoplus_{m\geq 0} (f_{2})_*\mathcal O_{X_2}(mD_2)$ 
is a finitely generated $\mathcal O_{\overline Y}$-algebra. 
We put 
$\overline X
:=\Proj _{\overline Y}\bigoplus _{m\geq 0}(f_{2})_*\mathcal 
O_{X_2}(mD_2)$. 
Then, $\overline f:\overline X\longrightarrow \overline Y$ has 
the required properties (i) and 
(ii) since $D_2$ is $f_2$-ample over $Y$. 
When $X$ is 
$\mathbb Q$-factorial, we replace $\overline X$ by 
its small projective $\mathbb Q$-factorialization 
(see \cite[Corollary 5.9]{fujino}). 
So, (1) holds. 
For (2), we apply Proposition \ref{propC} below. 
\end{proof}

The following is the blow-up whose exceptional divisor 
is the prescribed one. 

\begin{lem}\label{lemB}
Let $g:Z\longrightarrow X$ be a projective birational 
toric morphism. Let $E$ be an irreducible $g$-exceptional divisor 
on $Z$. We put
$$
h:X':=\Proj_X\bigoplus _{m\geq 0}g_*\mathcal O_Z(-mE)
\longrightarrow X
$$ 
and let $E'$ be the strict transform of $E$ on $X'$. 
Then $-E'$ is $h$-ample. 
So,  $X'\setminus E'\simeq X\setminus h(E')$.  

Furthermore, if $X$ is $\mathbb Q$-Gorenstein, that is, 
$K_X$ is $\mathbb Q$-Cartier, then 
$$K_{X'}=h^*K_{X}+aE',$$
where $a=a(E', X, 0)\in \mathbb Q$ is the discrepancy of $E'$ with 
respect to $(X,0)$ $($cf.~\cite[Definition 2.5]{km} and 
\cite[Definition 2.9]{part1}$)$.  
\end{lem}
\begin{proof}[Sketch of the proof]
Run the MMP (see \cite[3.1]{part1} or 
\ref{mmp2} below) over $X$ with respect to $-E$. 
In the notation of \ref{mmp2} below, $X'$ is the 
$(-E)$-canonical model over $X$. 
\end{proof}

\begin{prop}\label{propC}
Let $X$ be a toric variety and $\overline X$ 
an equivariant completion of 
$X$. Assume that $X$ has only terminal $($resp.
~canonical$)$ 
singularities. 
Then there exists a projective toric 
morphism $g:Z\longrightarrow \overline X$ such 
that $Z$ has only terminal $($resp.~canonical$)$ singularities and 
$g$ is isomorphic over $X$. 
Moreover, if $X$ is $\mathbb Q$-factorial, then 
we can make $Z$ to be $\mathbb Q$-factorial. 
\end{prop}
\begin{proof}
Let $h:V\longrightarrow \overline X$ be a projective toric resolution. 
We put $g:Z:=\Proj_{\overline X}
\bigoplus _{m\geq 0}h_*\mathcal O_V(mK_V)\longrightarrow 
\overline X$. Then $Z$ has only canonical singularities and 
$K_Z$ is $g$-ample. 
We note that $g$ is isomorphic over $X$. 
So, this $Z$ is a required one when $X$ has only canonical 
singularities. Thus, we may assume that $X$ has 
only terminal singularities. 
Since the number of the divisors that are exceptional over $Z$ 
and whose discrepancies are zero is finite, we can make 
$Z$ to have only terminal singularities by applying Lemma 
\ref{lemB} finitely many times. 

Furthermore, 
if $X$ is $\mathbb Q$-factorial, then we can make 
$Z$ to be $\mathbb Q$-factorial by \cite[Corollary 5.9]{fujino}. 
\end{proof}

The next proposition is useful when we 
treat non-$\mathbb Q$-factorial toric varieties. 

\begin{prop}\label{q-car}
Let $X$ be a toric variety and $D$ a Weil
divisor on $X$. 
Then there exists a small projective toric morphism 
$g:Z\longrightarrow X$ such that 
the strict transform $D_Z$ 
of $D$ on $Z$ is $\mathbb Q$-Cartier. 

Furthermore, let $U$ be the Zariski open set of $X$ 
on which $D$ is $\mathbb Q$-Cartier. 
Then we can construct $g:Z\longrightarrow X$ 
so that $D_Z$ is $g$-ample and 
$g$ is isomorphic over 
$U$.  
\end{prop}
\begin{proof}
By Corollary 5.8 in \cite{part1}, $\bigoplus _{m\geq 0}\mathcal 
O_X(mD)$ is a finitely generated $\mathcal O_X$-algebra. 
We put $g:Z:=\Proj_X\bigoplus _{m\geq 0}\mathcal O_X(mD)
\longrightarrow X$. This $g:Z\longrightarrow X$ has the 
required property. See, for example, \cite[Lemma 6.2]{km} 
or \cite[4.2 Proposition]{FA}. 
\end{proof}

\begin{cor}
Let $X$ be a toric variety. We assume that 
$X$ is $\mathbb Q$-Gorenstein, 
that is, $K_X$ is $\mathbb Q$-Cartier. 
Then there exists an equivariant completion $\overline X$ 
of $X$ such that $\overline X$ is $\mathbb Q$-Gorenstein. 
\end{cor}
\begin{proof}
Let $X'$ be an equivariant completion of $X$. 
We put $\overline X:=\Proj _{X'}\bigoplus _{m\geq 0}
\mathcal O_{X'}(mK_{X'})$. 
This $\overline{X}$ 
has the required property by Proposition \ref{q-car}. 
\end{proof}

The following theorem is a generalization of 
the elementary transformations. 
We need it for the MMP in \ref{mmp2} below. 

\begin{thm}[cf.~{\cite[Theorem 4.8]{part1}}]\label{g-flip}
Let $\varphi:
X\longrightarrow W$ be a projective birational toric morphism 
and $D$ a $\mathbb Q$-Cartier Weil divisor on $X$ such 
that $-D$ is $\varphi$-ample. 
We put 
$$\varphi^+:
X^+:=\Proj_W\bigoplus _{m\geq 0}\varphi_
*\mathcal O_X(mD)\longrightarrow W
$$ 
and let $D^+$ be the strict transform of $D_W$ on 
$X^+$, where $D_W:=\varphi_*D$. 
Then $\varphi^+$ 
is a small projective toric morphism such that 
$D^+$ is a $\varphi^+$-ample $\mathbb Q$-Cartier 
Weil divisor on $X^+$. 

Let $U$ be the Zariski open set of $W$ 
over which $\varphi$ is isomorphic. 
Then, so is $\varphi^+$ over $U$.  

The commutative diagram
$$
\begin{matrix}
X & \dashrightarrow & \ X^+ \\
{\ \ \ \ \searrow} & \ &  {\swarrow}\ \ \ \ \\
 \ & W &  
\end{matrix}
$$ 
is called the {\em{elementary transformation (with respect to 
$D$)}} if $\varphi:X\longrightarrow W$ is 
small $($cf.~\cite[Theorem 4.8]{part1}$)$. 
\end{thm}
\begin{proof}
We put $\varphi':X':
=\Proj _W\bigoplus _{m\geq 0}\mathcal O_W(mD_W)\longrightarrow W$ 
and let $D'$ be the strict transform of $D_W$ on $X'$ as in 
Proposition \ref{q-car}. 
Then, by the negativity lemma 
(see Lemma \ref{neg-lem} below), 
$\varphi'_*\mathcal O_{X'}(mD')
\simeq \varphi_*\mathcal O_X(mD)$ for every $m\geq 0$. 
We note that $\varphi'$ is small. 
Thus, we obtain 
\begin{eqnarray*}
X'&=&\Proj _W\bigoplus _{m\geq 0}\mathcal O_W(mD_W)\\
&\simeq &\Proj _W\bigoplus
_{m\geq 0}\varphi'_*\mathcal O_{X'}
(mD')\\
&\simeq &\Proj _W\bigoplus _{m\geq 0}\varphi
_*\mathcal O_{X}
(mD)=X^+. 
\end{eqnarray*} 
So, $\varphi^+$ and $D^+$ have the required properties. 
Note that this $X^+$ is the $D$-canonical model over $W$ 
in the notation of \ref{mmp2} below. 
See also Example \ref{4.3}.  
\end{proof}

Let us recall the following well-known negativity lemma (\cite[Lemma 
4.10]{part1}), 
which we already used in the proof of Theorem \ref{g-flip}. 
The proof can be found in \cite[Lemma 3.38]{km}. 

\begin{lem}[the Negativity Lemma]\label{neg-lem}
We consider a commutative diagram 
$$
\begin{matrix} 
&Z&\\
\ \ \ \ \swarrow & & \searrow \ \ \ \ \\
U & \dashrightarrow & \ V \\
{\ \ \ \ \searrow} & \ &  {\swarrow}\ \ \ \ \\
 \ & W &  
\end{matrix}
$$
and $\mathbb Q$-Cartier divisors $D$ and $D'$ 
on $U$ and $V$, respectively, where 
\begin{itemize}
\item[(1)] $f:U\longrightarrow W$ and 
$g:V\longrightarrow W$ are 
proper birational morphisms between normal varieties, 
\item[(2)] $f_*D=g_*D'$, 
\item[(3)] $-D$ is $f$-ample and 
$D'$ is $g$-ample, 
\item[(4)] $\mu:Z\longrightarrow U$, $\nu:Z\longrightarrow V$ are 
common resolutions. 
\end{itemize}
Then $\mu^*D=\nu^*D'+E$, 
where $E$ is an effective $\mathbb Q$-divisor and is exceptional 
over $W$. Moreover, 
if $f$ or $g$ is non-trivial, then $E\ne 0$.  
\end{lem}

\begin{rem}\label{rem1}
In Theorem \ref{g-flip}, 
let us further assume 
that $X$ is $\mathbb Q$-factorial and 
$\rho(X/W)=1$. 
If $\varphi$ 
contracts a divisor, then $W$ is $\mathbb Q$-factorial. 
In particular, $D_W$ is $\mathbb Q$-Cartier. 
So, $\varphi^+:X^+\longrightarrow W$ is an isomorphism. 
If $\varphi$ is small, then $X^+$ is $\mathbb Q$-factorial 
and $\rho(X^+/W)=1$. For the 
non-$\mathbb Q$-factorial case, see the examples in 
Section \ref{sa}.  
\end{rem}

The following Minimal Model 
Program (MMP, for short) for 
toric varieties is a slight generalization of the MMP 
explained in \cite[3.1]{part1}. 
This MMP works without the 
$\mathbb Q$-factoriality assumption. 
See also Remark \ref{itti} below. 

\begin{say}[Minimal Model Program for Toric Varieties]\label{mmp2}
We start with a projective 
toric morphism $f:X\longrightarrow Y$ and 
a $\mathbb Q$-Cartier divisor $D$ on $X$. 
Let $l$ be a positive integer such that $lD$ is a Weil divisor. 
We put $X_0:=X$ and $D_0:=D$. 
The aim is to set up a recursive procedure which 
creates intermediate $f_i:X_i\longrightarrow Y$ and $D_i$ 
on $X_i$. 
After finitely many steps, we obtain a finial objects 
$\widetilde {f}:\widetilde {X}\longrightarrow Y$ 
and $\widetilde {D}$. 
Assume that we already constructed $f_i:X_i\longrightarrow Y$ 
and $D_i$ with 
the following properties: 
\begin{enumerate}
\item[(i)] $f_i$ is projective, 
\item[(ii)] $D_i$ is a $\mathbb Q$-Cartier divisor on $X_i$. 
\end{enumerate}

If $D_i$ is $f_i$-nef, then we set $\widetilde{X}:=X_i$ 
and $\widetilde {D}:=D_i$. 
Assume that $D_i$ is not $f_i$-nef. 
Then we can take an 
extremal ray $R$ of $NE(X_i/Y)$ such that $R\cdot D_i<0$. 
Thus we have a contraction morphism $\varphi_R:X_i
\longrightarrow W_i$ over 
$Y$. 
If $\dim W_i <\dim X_i$, then we set $\widetilde{X}
:=X_i$ and $\widetilde{D}:=D_i$ 
and stop the process. 
If $\varphi_R$ is birational, 
then we put 
$$X_{i+1}:=\Proj_{W_i}\bigoplus_{m\geq 0}
\varphi_{R*}\mathcal O_{X_i}(mlD_i)
$$ and let 
$D_{i+1}$ be the strict transform of $\varphi_{R*}D_i$ on 
$X_{i+1}$ (see Theorem \ref{g-flip}). 
By counting the number of 
the torus invariant irreducible divisors, 
we may assume that $\varphi_R:X_i\longrightarrow W_i$ 
is small or $\dim W_i<\dim X_i$ after finitely many steps. 
By Theorem 4.9 (Termination of Elementary Transformations) in \cite{part1}, 
there are no infinite sequences of the elementary transformations with 
respect to $D_i$ (cf.~Theorem \ref{g-flip}). Therefore, 
this process always terminates and we 
obtain $\widetilde{f}:\widetilde{X}\longrightarrow Y$ 
and $\widetilde{D}$. 
We note that the relative Picard number may increase 
in the process (see Example \ref{fli} below). 
When $\widetilde D$ is $\widetilde f$-nef, $\widetilde X$ is 
called a {\em{$D$-minimal model over $Y$}}. 
We call this process ($D$-){\em{Minimal Model Program 
over}} $Y$, 
where $D$ is the divisor used in the process. 
When we apply the Minimal Model Program (MMP, for short), 
we say that, 
for example, we 
{\em{run the MMP over $Y$ with respect to the 
divisor $D$}}. 
If $\widetilde X$ is a $D$-minimal model over $Y$, 
then we put 
$$X^{\dagger}:=\Proj_Y \bigoplus_{m\geq 0}
{\widetilde f}_*\mathcal O_{\widetilde X}(ml\widetilde D).
$$ 
It is not difficult to see 
that $X^{\dagger}\simeq \Proj_Y \bigoplus_{m\geq 0}
{f}_*\mathcal O_{X}(mlD)$. We call $X^{\dagger}$ the 
{\em{$D$-canonical model over $Y$}}. 
We note that there exists a toric morphism 
$\widetilde X\longrightarrow X^\dagger$ over $Y$ 
which corresponds to 
${\widetilde f}^*{\widetilde f}_*\mathcal O_{\widetilde X}
(k\widetilde D)\longrightarrow \mathcal O_{\widetilde X}
(k\widetilde D)\longrightarrow 0$, 
where $k$ is a sufficiently large and divisible integer 
(see \cite[Proposition 4.1]{part1}).
\end{say}

\begin{rem}\label{itti}
(i) When $X$ 
is $\mathbb Q$-factorial, this process coincides with the 
one explained in \cite[3.1]{part1}. 
See Remark \ref{rem1}. 

(ii) If $X$ has only terminal (resp.~canonical) 
singularities and 
$D=K_X$, then so does $X_i$ 
for every $i$. It is an easy consequence of the negativity 
lemma (see Lemma \ref{neg-lem}). 
\end{rem}

The following Theorems \ref{kore} and 
\ref{supple-ne} are 
the main results in this paper. 
We divide them since Theorem \ref{kore} is sufficient 
for various applications and the proof of 
Theorem \ref{supple-ne} is complicated. 

\begin{thm}[Equivariant 
completions of toric contraction morphisms]\label{kore}
Let $f:X\longrightarrow Y$ be a projective toric morphism. 
Let $\varphi:=\varphi_R:X\longrightarrow W$ be 
the contraction morphism 
over $Y$ with respect to an 
extremal ray $R$ of $NE(X/Y)$. 
Then there exists an equivariant completion 
of $\varphi:X\longrightarrow W$ as follows{\em{;}} 
$$
\begin{matrix}
\overline \varphi: & \overline X & \longrightarrow &\overline W\\
                   & \cup & &\cup\\
\varphi: & X& \longrightarrow &W,\\
\end{matrix}
$$ 
where 
\begin{enumerate}
\item[(i)] $\overline X$ and $\overline W$ 
are equivariant completions of $X$ and $W$, and 
\item[(ii)] $\overline \varphi$ is a projective toric morphism with 
the relative Picard number $\rho (\overline X/\overline W)=1$. 
\end{enumerate} 
Furthermore, 
\begin{enumerate}
\item[(1)] if $X$ is $\mathbb Q$-factorial, then we can make 
$\overline X$ to be $\mathbb Q$-factorial, and 
\item[(2)] if $X$ has only $($$\mathbb Q$-factorial$)$ 
terminal $($resp.~canonical$)$ 
singularities and $-K_X$ is 
$\varphi$-ample, then we can make $\overline X$ to have 
only $($$\mathbb Q$-factorial$)$ 
terminal $($resp.~canonical$)$ singularities. 
\end{enumerate}

Let $\overline Y$ 
be an equivariant completion of $Y$. Then we 
can construct $\overline \varphi$ with the following 
property{\em{:}} 
\begin{enumerate}
\item[(iii)] $\overline W\longrightarrow \overline Y$ is an 
equivariant completion of $W\longrightarrow Y$ such that 
$\overline W\longrightarrow \overline Y$ is projective. 
\end{enumerate}
\end{thm}

\begin{proof} 
Let $W'$ be an equivariant completion of $W$. 
If $\overline Y$ is given, then we can take $W'$ to be 
projective over $\overline Y$ by Proposition \ref{sai}. 
Let $\varphi':X'\longrightarrow W'$ be an equivariant 
completion of $\varphi:X\longrightarrow W$. 
By Proposition 
\ref{sai}, we may assume that $\varphi'$ is projective. 
We may further assume that $X'$ is $\mathbb Q$-factorial 
(resp.~$X'$ has only 
($\mathbb Q$-factorial) terminal or canonical 
singularities) when $X$ is 
$\mathbb Q$-factorial (resp.~$X$ has only 
($\mathbb Q$-factorial) terminal or canonical singularities). 
Let $D$ be a $\mathbb Q$-Cartier divisor on $X$ such 
that $-D$ is $\varphi$-ample. 
Take a $\mathbb Q$-Cartier divisor $D'$ on $X'$ such 
that $D'|_X=D$. 
We note that we can always take such 
$D'$ by Proposition \ref{q-car} if we modify $X'$ 
suitably. 
We put $D'=K_{X'}$ in the case (2). 
Run the MMP (as explained in \ref{mmp2}) over $W'$ 
with respect to $D'$. 
If an extremal ray $R$ does not contain the numerical equivalence 
class of the 
curves contracted by $\varphi:X\longrightarrow W$, 
then the contraction with respect to $R$ occurs outside $X$. 
So, we obtain 
$$
X'=:X'_0\dashrightarrow X'_1\dashrightarrow X'_2
\dashrightarrow \cdots \dashrightarrow X'_k=:\overline X
$$ 
over $W'$ and a contraction $\overline \varphi:\overline X
\longrightarrow \overline W$ such that 
$\rho(\overline X/\overline W)=1$ and 
$\overline \varphi$ contracts the 
curves in the fibers of $\varphi$. 
It is easy 
to see that $\overline \varphi:\overline X\longrightarrow 
\overline W$ has the required properties. 
See also Remarks \ref{rem1} and \ref{itti}. 
\end{proof}

\begin{thm}\label{supple-ne}
We use the same notation as in {\em{Theorem \ref{kore}}}. 
We can generalize {\em{Theorem \ref{kore} (2)}} as follows$:$ 
\begin{enumerate}
\item[$(2')$] if $X$ has only $($$\mathbb Q$-factorial$)$ 
terminal $($resp.~canonical$)$ 
singularities and $-K_X$ is 
{\em{$\varphi$-nef}}, then we can make $\overline X$ to have 
only $($$\mathbb Q$-factorial$)$ 
terminal $($resp.~canonical$)$ singularities. 
\end{enumerate}
\end{thm}
\begin{proof}
By Theorem \ref{kore} (2), we may assume that 
$-K_X$ is not $\varphi$-ample, or equivalently, 
$K_X$ is $\varphi$-numerically trivial. 
As in the proof of Theorem \ref{kore}, we run the 
MMP over $W'$ with 
respect to $D'=K_{X'}$. 
In this case, 
we obtain 
$$
X'=:X'_0\dashrightarrow X'_1\dashrightarrow X'_2
\dashrightarrow \cdots \dashrightarrow X'_k=:\widetilde X
$$ 
over $W'$, and $\widetilde X$ is a 
$D'$-minimal model over $W'$, that is, $K_{\widetilde X}$ 
is nef over $W'$. 
It is easy to see that each step occurs outside $X$. 
Note that $K_{\widetilde X}$ is not ample over 
$W'$ since $K_X$ is $\varphi$-numerically 
trivial. 

Let $B$ be the complement of the big tours in $X$ regarded 
as a reduced divisor. 
Then it is well-known that 
$K_X+B\sim 0$. 
So, $B$ is $\varphi$-numerically trivial. 
Therefore, it is not difficult to see that there 
exists an effective torus-invariant Cartier divisor $E$ on 
$X$ such that $-E$ is $\varphi$-ample. 
Let $F$ be the closure of $E$ on $\widetilde X$. 
By modifying $\widetilde X$ birationally outside $X$ 
(if necessary), 
we may assume that $F$ is $\mathbb Q$-Cartier (see 
Proposition \ref{q-car}). 
Run the MMP over $W'$ with respect to $F$. 
For each step, we choose a $K$-trivial extremal ray $R$, 
that is, $K\cdot R=0$, where $K$ is the canonical divisor. 
Then we obtain a sequence 
$$
\widetilde X=:\widetilde X_0
\dashrightarrow \widetilde X_1\dashrightarrow 
\widetilde X_2
\dashrightarrow \cdots \dashrightarrow 
\widetilde X_l=:\overline X
$$ 
over $W'$ and a contraction $\overline \varphi:\overline 
X\longrightarrow \overline W$ 
such that $\overline \varphi$ contracts the curves in the 
fibers of $\varphi$. 
We note that $(\widetilde X, \varepsilon F)$ has only 
terminal singularities for $0\leq \varepsilon \ll 1$ 
(resp.~$\widetilde X$ has only canonical singularities) 
when $X'$ has only terminal (resp.~canonical) singularities. 
So, the pair 
$(\overline X, \varepsilon \overline F)$, 
where $\overline F$ is the strict transform of $F$, 
has only 
terminal singularities for $0\leq \varepsilon \ll 1$ 
(resp.~$\overline X$ has only canonical singularities) 
by \cite[Lemma 3.38]{km}. 
We note that each step of the above MMP 
does not contract any components 
of $F$ since it occurs outside $X$. 
Therefore, $\overline \varphi: \overline X\longrightarrow 
\overline W$ has the desired properties. 
\end{proof}

\begin{rem}
The assumptions on $K_X$ in Theorems \ref{kore} and 
\ref{supple-ne} are 
useful when we construct global (toric) examples 
of {\em{flips}} and {\em{flops}}.  
\end{rem}

The following is a question of J.~Koll\'ar. 

\begin{que}
Let $f:X\longrightarrow Y$ be a projective equivariant 
morphism between toric varieties 
with connected fibers. Assume that 
$\rho(X/Y)=k\geq 2$. 
Is it possible to compactify $f$ equivariantly preserving 
$\rho=k$?  
\end{que}

\section{Applications of equivariant completions}\label{sec-app}

In this section, we treat some applications of Theorem 
\ref{kore} and related topics. 

\begin{say} 
The next theorem is a direct consequence of 
Theorem \ref{kore} and Reid's description of the toric 
contraction 
morphisms. 
Theorem \ref{desc} was obtained by Reid when 
$X$ is {\em{complete}}. 
For the details, see, for instance, \cite[Corollary 14-2-2]{ma}, 
where Matsuki corrected minor errors in \cite{reid}. 
See \cite[Remark 14-2-3]{ma}. 
Remark \ref{321} below is a supplement to \cite[Corollary 14-2-2]{ma}. 
For the combinatorial aspects of this theorem, 
see \cite{sato}. 

\begin{thm}[cf.~{\cite[(2.5) Corollary]{reid}}]\label{desc}
Let $f:X\longrightarrow Y$ be a projective toric morphism. 
Assume that $X$ is $\mathbb Q$-factorial. 
Let $R$ be an extremal ray of $NE(X/Y)$ and 
$\varphi_R:X\longrightarrow W$ the contraction morphism 
over 
$Y$ with respect to $R$. 
Let 
$$
\begin{matrix}
 & A & \longrightarrow &B \\
                   & \cap & &\cap\\
\varphi_R: & X& \longrightarrow &W\\
\end{matrix}
$$ 
be the loci on which $\varphi_R$ is not an isomorphism{\em{;}} 
$A$ and $B$ are irreducible, 
$\varphi^{-1}(P)_{\red}$ is a $\mathbb Q$-factorial 
projective toric $(\dim A-\dim B)$-fold with 
the Picard number one for every point $P\in B$.  
More precisely, 
there exist an open covering $B=\bigcup _{i\in I}U_i$ and 
a $\mathbb Q$-factorial projective 
toric variety $F$ with 
the Picard number $\rho (F)=1$ such that 
\begin{enumerate}
\item[(i)] $U_i$ is a torus invariant open 
subvariety of $B$ for 
every $i$, 
\item[(ii)] there exists a finite toric morphism 
$U'_i\longrightarrow U_i$ such that 
$$
(U'_i\times _B A)^{\nu}\simeq U'_i\times F
$$
for every $i\in I$, where $(U'_i\times _B A)^{\nu}$ 
is the normalization of $U'_i\times _B A$. 
\end{enumerate} 
We note that 
$-K_F$ is an ample $\mathbb Q$-Cartier divisor since 
$\rho(F)=1$. 
\end{thm}
\begin{proof} 
By Theorem \ref{kore}, we obtain an equivariant 
completion:  
$$
\begin{matrix}
\overline \varphi: & \overline X & \longrightarrow &\overline W\\
                   & \cup & &\cup\\
\varphi_R: & X& \longrightarrow &W.\\
\end{matrix}
$$ 
We may assume that $\overline X$ is $\mathbb Q$-factorial, 
$\overline\varphi$ is projective, and 
$\rho(\overline X/\overline W)=1$. 
Let 
$$
\begin{matrix}
 & \overline A & \longrightarrow &\overline B \\
                   & \cap & &\cap\\
\overline\varphi: & \overline X& \longrightarrow &\overline W\\
\end{matrix}
$$ 
be the loci on which $\overline \varphi$ is not 
an isomorphism. 
Apply Reid's description:~\cite[(2.5) Corollary]{reid} 
to $\overline\varphi$. For the detailed description of 
$\overline{\varphi}:\overline{X}\longrightarrow \overline{W}$, 
see \cite[Corollary 14-2-2]{ma} and Remark \ref{321} below. 
\end{proof}

In the above theorem, the assumption that 
$X$ is $\mathbb Q$-factorial plays a crucial role. 
See Example \ref{sa-ex} below and \cite[Example 1.1]{fs-flop2}. 

\begin{rem}[Supplements to the 
description of contractions of extremal rays by Matsuki]\label{321}
In this remark, we use the same notation as in 
\cite[Chapter 14]{ma}. 
In \cite[Corollary 14-2-2]{ma}, Matsuki claimed that 
$E\times _F U'_{\tau(w)_Y}\cong G\times U'_{\tau(w)_Y}$. 
In our notation in Theorem \ref{desc}, he claims that 
$U'_i\times _B A\simeq U'_i\times F$. 
However, it is not true in general. 
We have to take the normalization of the left hand side. So, the correct 
statement should be 
$(E\times _F U'_{\tau(w)_Y})^{\nu}\cong G\times U'_{\tau(w)_Y}$, where 
$(E\times _F U'_{\tau(w)_Y})^{\nu}$ is the normalization of 
$E\times _F U'_{\tau(w)_Y}$. We note that 
$(E\times _F U'_{\tau(w)_Y})^{\nu}$ is irreducible since 
$E\longrightarrow F$ has connected fibers. 
For the details, see \cite[Lemma 5.6]{karu}. 
Therefore, $\varphi^{-1}_{R}(P)_{\red}$ is not necessarily isomorphic to 
$G$. 
Let $O(\gamma)\subset F$ 
be the orbit associated to a cone $\gamma$. 
Then $\varphi^{-1}_{R}(O(\gamma))_{\red}\simeq 
G_\gamma\times O(\gamma)$ and $\varphi^{-1}_{R}
(O(\gamma))_{\red}\longrightarrow O(\gamma)$ is isomorphic 
to the second projection $G_\gamma\times O(\gamma)\longrightarrow 
O(\gamma)$, 
where $G_\gamma$ is an $(n-\beta)$-dimesnional $\mathbb Q$-factorial 
projective toric variety with the Picard number 
$\rho(G_\gamma)=1$, 
Note that $G_\gamma$ is 
defined by $n-\beta+1$ one-dimensional vectors $\{v_{\beta+1}, 
\cdots, v_{n+1}\}$ for any $\gamma$. 
However, in general, $G_{\gamma_1}\not\simeq G_{\gamma_2}$ 
for two distinct cones $\gamma_1, \gamma_2$. 
It is because the lattice group that defines $G_\gamma$ depends on $\gamma$. 
So, $E\longrightarrow F$ is not necessarily 
a fiber bundle but a {\em{quasi}}-fiber 
bundle in Ishida's notation.  
The following example may help the reader to understand it. 
\end{rem}
\begin{ex}[Extremal Fano contraction] 
We fix $N=\mathbb Z^3$ and $N'=\mathbb Z$.
We put
\begin{align*}
v_1  &= (0,0,1), & v_2&=(-1,0,0), & v_3&=(1,0,-1),\\
v_4  &= (0,-1,0), & v_5&=(0,2,-1). &
\end{align*} 
We consider the following fan. 
$$
\Delta=
\left \{
\begin{array}{cccc}
\langle v_1, v_2, v_4\rangle, &
\langle v_1, v_2, v_5\rangle, &
\langle v_1, v_3, v_4\rangle, &
\langle v_1, v_3, v_5\rangle, \\
\langle v_2, v_3, v_4\rangle, &
\langle v_2, v_3, v_5\rangle, & 
\text{and their faces} &
\end{array}
\right \}. $$ 
We define $X=X(\Delta)$. Then $X$ is a $\mathbb Q$-factorial
projective toric $3$-fold with $\rho(X)=2$.
Let $N_{\mathbb R}\simeq \mathbb R^3\longrightarrow N'_{\mathbb R}\simeq
\mathbb R$ be the projection to the second coordinate.
It induces a toric morphism $f:X\longrightarrow \mathbb P^1$.
Then $f:X\longrightarrow \mathbb P^1$ is an extremal contraction.
We have the following properties:
\begin{itemize}
\item[(i)] $f^{-1}(0)$ is non-reduced since $v_5$ is mapped onto $2\in N'$,
\item[(ii)] $f^{-1}(0)_{\red}$ is isomorphic to a weighted projective space
$\mathbb P (1,1,2)$, and
\item[(iii)]$g:=f|_{X\setminus f^{-1}(0)}:Y:=X\setminus
f^{-1}(0)\longrightarrow
Z:=\mathbb P^1\setminus \{0\}$ is isomorphic to the second projection
$\mathbb P^2
\times \mathbb A^1\longrightarrow \mathbb A^1$.
\end{itemize}
We note that $-K_X\cdot C\geq 3$ for any curve $C$ in the fibers
$g:Y\longrightarrow Z$.
On the other hand, there exists a torus invariant curve $C_0$ in the
fiber $f^{-1}(0)$ such
that $-K_X\cdot C_0=\frac{3}{2}$.
It can be checked by adjunction and the 
computations in \cite[Section 2]{fujino} 
(cf.~Theorem \ref{length}).
\end{ex}

\begin{rem}
In Theorem \ref{desc}, let $\Delta$ be the fan such that 
$X=X(\Delta)$. Then $\Delta$ need not 
contain $n$-dimensional cones, where $n=\dim X$. 
\end{rem}
\begin{rem}
Let $f:X\longrightarrow Y$ be a projective equivariant 
morphism between toric varieties. 
Let $\varphi_R:X\longrightarrow W$ be the extremal 
contraction associated to an extremal ray $R$ of $NE(X/Y)$. 
Assume that $X$ is $\mathbb Q$-factorial. 
If $X$ is complete, then Reid obtained the combinatorial 
descriptions of $\varphi_R$ in \cite{reid} by using the notion of 
{\em{walls}}. 
Sato generalized Reid's combinatorial 
descriptions for non-complete toric 
varieties in \cite{sato} by using Theorem \ref{kore} and 
the notion of {\em{extremal primitive relations}}. 
Examples in Section \ref{sa} and Example 1.1 in 
\cite{fs-flop2} 
imply that it is impossible to 
describe toric extremal contractions 
combinatorially without $\mathbb Q$-factoriality.  
\end{rem}

\begin{rem}\label{rigidity}
In \cite[Chapter 14]{ma}, Matsuki left the details 
of the verifications for the relative case to the reader 
in various places. 
In the relative case in \cite[Proposition 14-1-2]{ma}, 
all we need is the rigidity lemma (see, 
for example, \cite[Lemma 1.6]{km}). The rest are straightforward. 
In \cite[Chapter 14]{ma}, $X(\Delta)$ is always assumed to be 
{\em{complete}} even in the relative setting of $\phi:X(\Delta)\longrightarrow 
S(\Delta_S)$. So, there are no difficulties to handle the relative setting in \cite
[Chapter 14]{ma}. For the {\em{true}} relative setting, that is, 
$\phi:X(\Delta)\longrightarrow S(\Delta _S)$ is a projective morphism 
and 
$X(\Delta)$ is not necessarily complete, see \cite{part1}.  
\end{rem}

Here 
is a general remark on equivariant completions 
of toric varieties. 

\begin{rem}\label{ippan}
Let $X$ be a toric variety corresponding to a 
fan $\Delta$. It is well-known that compactifying $X$ 
equivariantly is equivalent to compactifying $\Delta$. 
We know that to compactify 
$\Delta$ without Sumihiro's theorem 
is very difficult. 
Recently, Ewald and Ishida independently succeeded in 
compactifying (not necessarily rational) fans without 
using Sumihiro's theorem (see \cite{ei}). 
\end{rem}
\end{say}
\begin{say}
In \cite[Corollary 4.6]{app}, we proved that the 
target space of a Mori fiber space has at most log 
terminal singularities. 
In dimension three, it is conjectured 
that the target space has only canonical singularities 
(see, for instance, \cite[Conjecture 0.2]
{pro}). 
Before we explain an example of Mori fiber spaces, 
let us recall the definition of the Mori fiber space. 
\begin{defn}[Mori fiber space]
A normal projective variety $X$ with only $\mathbb Q$-factorial 
terminal singularities with a morphism $\Phi:X\longrightarrow Y$ 
is a {\em{Mori fiber space}} if (i) $\Phi$ is a morphism 
with connected fibers onto a normal projective variety $Y$ of $\dim Y<\dim X$, 
(ii) $-K_X$ is $\Phi$-ample, and (iii) $\rho(X/Y)=1$. 
\end{defn}
The following is an example of $4$-dimensional 
Mori fiber spaces. 

\begin{ex}[Mori fiber space whose 
target space has a bad singularity]\label{morifiber}
Let $\mathbb Z_4=\langle \zeta\rangle$ be the 
cyclic group of fourth roots of unity 
with $\zeta=\sqrt{-1}$. 
Let $\mathbb P^1\times \mathbb C^3\longrightarrow 
\mathbb C^3$ be the second projection. 
We consider the following actions of $\mathbb Z_4$ on 
$\mathbb P^1\times \mathbb C^3$ and 
$\mathbb C^3$: 
\begin{eqnarray*}
\left( [u:v], (x,y,z)\right)
&\longrightarrow &\left( [u: -v], (\zeta x, 
\zeta y, \zeta z) \right),\\
(x,y,z)&\longrightarrow &(\zeta x, 
\zeta y, \zeta z), 
\end{eqnarray*}
where $[u:v]$ is the homogeneous coordinate of 
$\mathbb P^1$. 
We put $X:=(\mathbb P^1\times \mathbb C^3)\slash 
\mathbb Z_4$ and $Y:=\mathbb C^3\slash 
\mathbb Z_4$. 
Then the induced equivariant morphism 
$f:X\longrightarrow Y$ has the 
following properties: 
\begin{enumerate}
\item[(i)] $X$ has terminal quotient singularities along the 
central fiber of $f$,  
\item[(ii)] $Y$ has a $\frac{1}{4}(1,1,1)$ quotient 
singularity, which is not canonical, 
\item[(iii)] $X$ and $Y$ are $\mathbb Q$-factorial, 
\item[(iv)] $\rho(X/Y)=1$, and 
\item[(v)] $-K_X$ is $f$-ample. 
\end{enumerate} 
By applying Theorem \ref{kore} (2), we obtain a toric 
Mori fiber space $\overline f:\overline{X}\longrightarrow 
\overline{Y}$ that is an equivariant completion of 
$f:X\longrightarrow Y$. 
Note that we can make $\overline Y$ projective 
by Theorem \ref{kore} (iii). Thus, $\overline f: 
\overline X\longrightarrow \overline Y$ is a Mori 
fiber space such that the target space $\overline Y$ 
has a singularity that is not canonical.  
\end{ex}

This example shows that our theorem is useful when 
we construct {\em{global}} examples from local ones. 
For a more combinatorial treatment, 
see \cite[Example 5.1]{fs-flop}. 
\end{say}
\begin{say} 
The final 
theorem is a slight generalization of \cite[Theorem 0.1]{fujino}. 

\begin{thm}[Length of an extremal ray]\label{length}  
Let $f:X\longrightarrow Y$ be a projective surjective 
equivariant morphism between toric varieties. 
Let $D=\sum _jd_jD_j$ be a $\mathbb Q$-divisor, 
where $D_j$ is an irreducible torus invariant divisor and 
$0\leq d_j\leq 1$ for every $j$. 
Assume that $K_X+D$ is $\mathbb Q$-Cartier. 
Then, for each extremal ray $R$ of $NE(X/Y)$, 
there exists an irreducible curve $C$ such that 
$[C]\in R$ and 
$$
-(K_X+D)\cdot C\leq \dim X+1. 
$$ 
Moreover, we can choose $C$ in such a way that 
$$
-(K_X+D)\cdot C\leq \dim X
$$
unless $X\simeq \mathbb P^{\dim X}$ and $\sum _jd_j<1$. 
Here, we do not claim that $C$ is a torus invariant 
curve. We note 
that $R$ may contain no numerical 
equivalence classes of torus invariant curves. 
\end{thm}
\begin{proof}[Sketch of the proof]
If $Y$ is a point, then this is the main 
theorem of \cite{fujino}. 
So, we may assume that $\dim Y\geq 1$. 
Since the arguments in Step 2 in the proof of 
\cite[Theorem 0.1]{fujino} work with minor modifications, 
we may further assume that $X$ is $\mathbb Q$-factorial. 
Let $R$ be a $(K_X+D)$-negative 
extremal ray of $NE(X/Y)$. 
We consider the contraction $\varphi_R:X
\longrightarrow W$ over $Y$ with respect to $R$. 
Let $U$ be a quasi-projective torus invariant open subvariety 
of $W$ such that $X_U:=\varphi_{R}^{-1}(U)\longrightarrow 
U$ is not an isomorphism. 
It is not difficult to see that $X_U$ is $\mathbb Q$-factorial and 
$\rho(X_U/U)=1$ (see \cite[Example 1.1]{fs-flop2}). 
We note that $\Pic(X)\otimes \mathbb Q\longrightarrow 
\Pic (X_U)\otimes \mathbb Q$ is surjective. 
So, by shrinking $W$, we may assume that $X$ and 
$W$ are quasi-projective. 
By Theorem \ref{kore}, 
we have an equivariant completion of 
$\varphi:=\varphi_R:X\longrightarrow W$, 
that is, 
$$
\begin{matrix}
\varphi: & X & \longrightarrow &W \\
                   & \cap & &\cap\\
\overline 
\varphi: & \overline X& \longrightarrow &\overline W,\\
\end{matrix}
$$ 
where $\overline X$ and $\overline W$ are 
$\mathbb Q$-factorial projective toric varieties 
and $\rho (\overline 
X/\overline W)=1$. 
Let $\overline D$ be the closure of $D$ on 
$\overline X$. 
Then $-(K_{\overline X}+\overline D)$ is 
$\overline \varphi$-ample. 
Therefore, $\overline \varphi$ is the contraction morphism 
with respect to a suitable $(K_{\overline X}+
\overline D)$-negative extremal ray $Q\subset NE(\overline X/
\overline W)\subset NE(\overline X)$ 
(see \cite[(1.5)]{reid} and \cite[3.8]{fs-flop}). So, 
we can apply the arguments in Step 1 in the 
proof of \cite[Theorem 0.1]{fujino} to 
$\overline \varphi:\overline X\longrightarrow \overline W$. 
Let
$$
\begin{matrix}
 & A & \longrightarrow &B \\
                   & \cap & &\cap\\
\varphi_R: & X& \longrightarrow &W\\
\end{matrix}
$$ 
be the loci on which $\varphi_R$ is not an isomorphism. 
Let $\overline A$ (resp.~$\overline B$) be the closure of $A$ (resp.~$B$) in $\overline X$ (resp.~$\overline W$).
We can calculate $K_{\overline A}$ by adjunction (cf.~the computation of
$K_P$
in \cite[Proof of The Theorem]{fujino}). 
Let $F$ be a general fiber of $\overline A\longrightarrow \overline B$. 
Then $F$ is a $\mathbb Q$-factorial projective toric $(\dim A-\dim B)$-fold with the
Picard number one
(cf.~Theorem \ref{desc}). 
So, it is sufficient to prove the following claim. 
\begin{claim}
There exists a curve $C$ in $F$ such that 
$$
-(K_X+D)\cdot C\leq -K_{\overline A}\cdot C=-K_F\cdot C\leq 
\dim A-\dim B+1\leq \dim X.  
$$
\end{claim}
If $C$ is in $F$, then the first inequality follows from the 
computations similar to the ones 
in Step 1 in \cite[Proof of The Theorem]{fujino}.
By adjunction, $K_{\overline A}|_{F}=K_F$. Thus, it is obvious that
$-K_{\overline A}\cdot
C=-K_{F}\cdot C$ for $C$ in $F$. 
The computaitons in \cite[Section 2]{fujino} imply the
existence of $C$ on $F$ 
such that $-K_F\cdot C\leq \dim F+1$.
\end{proof}
\end{say}

\section{Examples of non-$\mathbb 
Q$-factorial contractions}\label{sa}

In this section, we explain various examples of 
non-$\mathbb Q$-factorial 
toric contraction morphisms. 
All the examples are three-dimensional. 

The first one is a beautiful example due to Sato 
of divisorial contractions. 
This implies that it is difficult to describe the local 
behavior 
of divisorial contractions without 
the $\mathbb Q$-factoriality assumption 
even if the relative Picard number is one. 

\begin{ex}[Sato's non-$\mathbb Q$-factorial 
divisorial contraction]\label{sa-ex}
Let $e_1, e_2, e_3$ form the usual basis of 
$\mathbb Z^3$, and let $e_4$ be given by 
$$
e_1+e_2=e_3+e_4. 
$$ 
We put 
$$e_5=e_1+e_2=e_3+e_4$$ 
and 
$$e_6=e_2+e_3.$$ 
Let 
$$\Delta_Y=\{\langle e_1, e_2,e_3,e_4\rangle, 
\text{and its faces}\}$$ 
and $Y:=X(\Delta_Y)$. 
We put 
$$\Delta_X=\{\langle e_1, e_4, e_5\rangle, 
\langle e_1, e_3, e_5, e_6\rangle, 
\langle e_2, e_4, e_5, e_6\rangle, 
\text{and their faces}\}.$$ 
We define $X:=X(\Delta_X)$. 
Then $f:X\longrightarrow Y$ has the following properties: 
\begin{enumerate}
\item[(i)] $X$ has terminal singularities, 
\item[(ii)] $X$ is not $\mathbb Q$-factorial, 
\item[(iii)] $f$ is a projective birational equivariant morphism with 
$\rho(X/Y)=1$, 
\item[(iv)] $-K_X$ is $f$-ample, and 
\item[(v)] the exceptional locus contains a reducible divisor. 
\end{enumerate} 
Figure 1 helps us to understand the 
above contraction morphism. 
\newpage
\begin{figure}[hbtp]
\begin{center}
\setlength{\unitlength}{2.6mm}
\begin{picture}(40,42)
\put(16,4){\line(1,0) {6}}
\put(16,4){\line(0,1) {6}}
\put(22,4){\line(0,1) {6}}
\put(16,10){\line(1,0) {6}}

\put(30,18){\line(1,1) {6}}
\put(30,18){\line(0,1) {6}}
\put(30,18){\line(1,0) {6}}
\put(36,18){\line(0,1) {6}}
\put(30,24){\line(1,0) {6}}
\put(30,24){\line(1,-1) {6}}

\put(18,34){\line(1,1) {6}}
\put(18,34){\line(0,1) {6}}
\put(18,34){\line(1,0) {6}}
\put(24,34){\line(0,1) {6}}
\put(18,40){\line(1,0) {6}}
\put(18,40){\line(1,-1) {6}}
\put(21,37){\line(1,0){3}}

\put(4,28){\line(1,1) {3}}
\put(4,28){\line(0,1) {6}}
\put(4,28){\line(1,0) {6}}
\put(10,28){\line(0,1) {6}}
\put(4,34){\line(1,0) {6}}
\put(4,34){\line(1,-1) {6}}
\put(7,31){\line(1,0){3}}

\put(4,16){\line(1,1) {3}}
\put(4,16){\line(0,1) {6}}
\put(4,16){\line(1,0) {6}}
\put(10,16){\line(0,1) {6}}
\put(4,22){\line(1,0) {6}}
\put(4,22){\line(1,-1) {3}}
\put(7,19){\line(1,0){3}}

\put(29,17){\vector(-1,-1){5}}
\put(25,32){\vector(2,-3){4}}
\put(16,37){\vector(-2,-1){5}}
\put(7,27){\vector(0,-1){4}}
\put(8,14){\vector(1,-1){4}}

\put(18,2){$Y$}
\put(14,10){$e_1$}
\put(22,10){$e_3$}
\put(22,3){$e_2$}
\put(14,3){$e_4$}
\put(34,21){$e_5$}
\put(25,37){$e_6$}
\put(2,18){$X$}
\put(8,11){$f$}
\put(27, 13){$\varphi_1$}
\put(27,30){$\varphi_2$}
\put(32,16){$X_1$}
\put(20,32){$X_2$}
\put(1,30){$X_3$}
\put(8, 25){$\varphi_4$}
\put(13, 37){$\varphi_3$}
\end{picture}
\end{center}
\label{e2}
\caption{}
\end{figure}
We can easily check the following properties: 
\begin{enumerate}
\item[(1)] $X_1$ and $X_2$ are non-singular, 
\item[(2)] $\varphi_1$ and $\varphi_2$ are blow-ups, 
\item[(3)] $\varphi_3$ and $\varphi_4$ are flopping 
contractions, 
that is, $K_{X_i}$ is $\varphi_{i+1}$-numerically trivial for 
$i=2, 3$, 
\item[(4)] $X_3$ and $X$ are not $\mathbb Q$-factorial, 
\item[(5)] $X_3$ and $X$ have only terminal singularities. 
\item[(6)] $\rho(X_1/Y)=2$ and 
$\rho(X_2/X_1)=\rho(X_2/X_3)=\rho(X_3/X_4)=1$. 
\end{enumerate}
The ampleness of $-K_X$ follows from the 
convexity of the roofs of the 
maximal cones in $\Delta_X$ 
(cf.~\cite[(4.3) Proposition]{reid}). 
\end{ex}

The next is an example of {\em{flips}}. 
In this example, the relative Picard number increases 
by a flip. 

\begin{ex}[Non-$\mathbb Q$-factorial flip]\label{fli}
Let $e_1, e_2, e_3$ form the usual basis of 
$\mathbb Z^3$, and let $e_4$ be given by 
$$
e_1+e_2=e_3+e_4. 
$$ 
We put 
$f_1=(3,1,-2), f_2=(-1,1,2)\in \mathbb Z^3$. 
We consider the following fans: 
\begin{eqnarray*}
\Delta_a&=&\{\langle e_1, e_3, f_1, f_2\rangle, 
\langle e_2, e_4, f_1, f_2\rangle, 
\text{and their faces}\},\\ 
\Delta_b&=&\{\langle e_1, e_4, f_1\rangle, 
\langle e_2, e_3, f_2\rangle, 
\langle e_1, e_2, e_3, e_4\rangle,  
\text{and their faces}\}, {\text{and}}\\
\Delta_c&=&\{\langle e_1, e_2, e_3, e_4, f_1, f_2\rangle, 
\text{and its faces}\}. 
\end{eqnarray*}
We put $X:=X(\Delta_a), X^+:=X(\Delta_b)$, and 
$Y:=X(\Delta_c)$. 
Then we have a commutative diagram: 
$$
\begin{matrix}
X & \dashrightarrow & \ X^+ \\
{\ \ \ \ \ \searrow} & \ &  {\swarrow}\ \ \ \ \\
 \ & Y &  
\end{matrix}
$$
such that 
\begin{enumerate}
\item[(i)] $f:X\longrightarrow Y$ and 
$f^{+}:X^+\longrightarrow Y$ 
are both small projective equivariant morphisms, 
\item[(ii)] $\rho(X/Y)=1$ and $\rho(X^+/Y)=2$, 
\item[(iii)] $X$ and $X^+$ are not $\mathbb Q$-factorial, and 
\item[(iv)] $-K_X$ is $f$-ample and $K_{X^+}$ is $f^+$-ample. 
\end{enumerate}
Thus, this diagram is a so-called {\em{flip}}. 
Figure 2 helps us to understand this example. 
\begin{figure}[hbtp]
\begin{center}
\setlength{\unitlength}{2mm}
\begin{picture}(50,34)
\put(18,6){\line(1,1) {4}}
\put(18,6){\line(1,-1) {4}}
\put(22,2){\line(1,0) {8}}
\put(22,10){\line(1,0) {8}}
\put(34,6){\line(-1,1) {4}}
\put(34,6){\line(-1,-1) {4}}

\put(32,26){\line(1,1) {4}}
\put(32,26){\line(1,-1) {4}}
\put(36,30){\line(1,0) {8}}
\put(36,22){\line(1,0) {8}}
\put(44,22){\line(1,1) {4}}
\put(44,30){\line(1,-1) {4}}
\put(36,22){\line(0,1) {8}}
\put(44,22){\line(0,1) {8}}

\put(4,26){\line(1,1) {4}}
\put(4,26){\line(1,-1) {4}}
\put(8,30){\line(1,0) {8}}
\put(8,22){\line(1,0) {8}}
\put(16,22){\line(1,1) {4}}
\put(16,30){\line(1,-1) {4}}
\put(4,26){\line(1,0) {16}}

\put(16,20){\vector(1,-2){4}}
\put(36,20){\vector(-1,-2){4}}

\put(7,31){$e_1$}
\put(16,21){$e_2$}
\put(16,31){$e_3$}
\put(7,21){$e_4$}
\put(2,26){$f_1$}
\put(20,26){$f_2$}

\put(12,19){$X$}
\put(25,5){$Y$}
\put(39,25){$X^+$}
\put(16,15){$f$}
\put(35,15){$f^+$}
\end{picture}
\end{center}
\label{e3}
\caption{}
\end{figure}

We note that the points $e_1, e_3, f_1$, and 
$f_2$ are on $H_1=\{(x,y,z) | x+z=1\}$, 
$e_2, e_4, f_1$, and 
$f_2$ are on $H_2=\{(x,y,z) | y=1\}$. 
It is obvious that $H_3=\{(x,y,z) | x+y+z=2\}$ contains 
$f_1$ and $f_2$, while 
the points $e_1, e_2, e_3$, and 
$e_4$ are on $H_4=\{(x,y,z) | x+y+z=1\}$. 
We note that the non-trivial lattice points in the 
{\em{shed}} (see \cite[p.414 Definition]{reid}) of 
$\Delta_a$ are 
\begin{equation*}
\frac{1}{4} f_1+\frac{3}{4}f_2, 
\ \frac{1}{2}
f_1+\frac{1}{2}f_2, \ \text{and} \ \ 
\frac{3}{4} f_1+\frac{1}{4}f_2,\in \mathbb Z^3. 
\end{equation*}
We can check the following properties: 
\begin{enumerate}
\item[(1)] The flipping locus is $\mathbb P^1$ and $X$ has only 
canonical singularities (see \cite[(1.11) Definition]{reid}).  
\item[(2)] The flipping curve is contained in the singular locus of 
$X$.  
\item[(3)] $X^+$ has only one singular point, which is 
an ordinary double point. In particular, 
$X^+$ has only terminal singularities.  
\item[(4)] The flipped locus is $\mathbb P^1\cup \mathbb P^1$ 
and these two $\mathbb P^1$'s intersect each other 
at the singular point of $X^+$. 
\end{enumerate} 
The ampleness of $-K_X$ (resp.~$K_{X^+}$) 
follows from the convexity (resp.~concavity) 
of the roofs of the maximal 
cones (cf.~\cite[(4.3) Proposition]
{reid}). 
\end{ex}

The final example is a non-$\mathbb Q$-factorial divisorial 
contraction whose target space is not 
$\mathbb Q$-Gorenstein. 

\begin{ex}\label{4.3}
We use the same notation as in Example \ref{fli}. 
We put $f_3=(0,1,1)$. 
We note that 
$$
f_3=e_2+e_3=\frac{1}{4}f_1+\frac{3}{4}f_2. 
$$ 
We consider the following fans: 
\begin{eqnarray*}
\Delta_d&=&\{\langle e_1, e_4,f_3\rangle, 
\langle e_1, e_3, f_2, f_3\rangle, 
\langle e_2, e_4, f_2, f_3\rangle,  
\text{and their faces}\}, \\
\Delta_e&=&\{\langle e_1, e_2, e_3, e_4\rangle, 
\langle e_2, e_3, f_2\rangle, 
\text{and their faces}\}, {\text{and}}\\
\Delta_f&=&\{\langle e_1, e_2, e_3, e_4, f_2\rangle, 
\text{and its faces}\}. 
\end{eqnarray*}
We define $V:=X(\Delta _d)$, $V^+:=X(\Delta _e)$, and 
$W:=X(\Delta _f)$. 
Then we obtain a commutative diagram: 
$$
\begin{matrix}
V & \dashrightarrow & \ V^+ \\
{\ \ \ \ \ \searrow} & \ &  {\swarrow}\ \ \ \ \\
 \ & W &  
\end{matrix}
$$
such that 
\begin{enumerate}
\item[(i)] $\varphi:V\longrightarrow W$ is a projective birational 
equivariant morphism and $\varphi$ contracts a divisor, 
\item[(ii)] $\varphi^+:V^+\longrightarrow W$ is a small projective 
equivariant morphism, 
\item[(iii)] $-K_X$ is $\varphi$-ample and 
$K_{V^+}$ is $\varphi^+$-ample, 
\item[(iv)] $\rho(V/W)=\rho(V^+/W)=1$, 
\item[(v)] $V$ and $V^+$ have 
only terminal singularities, 
\item[(vi)] all $V$, $V^+$, and $W$ are not $\mathbb Q$-factorial, 
and 
\item[(vii)] $W$ is not $\mathbb Q$-Gorenstein. 
\end{enumerate}
See Figure 3. 
\begin{figure}[hbtp]
\begin{center}
\setlength{\unitlength}{2mm}
\begin{picture}(42,34)
\put(16,2){\line(1,0) {8}}
\put(16,2){\line(0,1) {8}}
\put(16,10){\line(1,0) {8}}
\put(28,6){\line(-1,-1) {4}}
\put(28,6){\line(-1,1) {4}}

\put(4,22){\line(1,0) {8}}
\put(4,22){\line(0,1) {8}}
\put(4,30){\line(1,0) {8}}
\put(4,22){\line(2,1) {8}}
\put(4,30){\line(2,-1) {8}}
\put(12,26){\line(1,0) {4}}
\put(12,22){\line(1,1) {4}}
\put(12,30){\line(1,-1) {4}}

\put(28,22){\line(0,1) {8}}
\put(28,22){\line(1,0) {8}}
\put(28,30){\line(1,0) {8}}
\put(36,30){\line(0,-1) {8}}
\put(40,26){\line(-1,1) {4}}
\put(40,26){\line(-1,-1) {4}}

\put(12,20){\vector(1,-2){4}}
\put(30,20){\vector(-1,-2){4}}

\put(3,31){$e_1$}
\put(12,21){$e_2$}
\put(12,31){$e_3$}
\put(3,21){$e_4$}
\put(11,27){$f_3$}
\put(16,27){$f_2$}

\put(1,25){$V$}
\put(20,5){$W$}
\put(31,25){$V^+$}
\put(12,15){$\varphi$}
\put(29,15){$\varphi^+$}
\end{picture}
\end{center}
\label{e4}
\caption{}
\end{figure}

We note that 
the small morphism $\varphi^+:V^+\longrightarrow W$ is 
the one given in Theorem \ref{g-flip}. 
This operation $V\dashrightarrow V^+$ 
preserves the relative Picard number over $W$. 
Note that the number of the torus invariant divisors decreases. 
\end{ex}

This example shows that we need to modify $W$ to continue the 
MMP even if $\varphi$ contracts a divisor (see \ref{mmp2}). 

In Examples \ref{sa-ex}, \ref{fli}, and \ref{4.3}, 
the 
varieties are not complete. 
To produce global examples, we just compactify them 
by Theorem \ref{kore}. 
More concrete global examples can be found in the 
appendix \cite{fs-flop2}. 
\ifx\undefined\bysame
\newcommand{\bysame|{leavemode\hbox to3em{\hrulefill}\,}
\fi

\end{document}